\documentclass[11pt]{amsart}
\usepackage{amscd,amssymb}
\usepackage{amsthm,amsmath,amssymb}
\usepackage[matrix,arrow]{xy}

\newtheorem{theorem}[equation]{Theorem}

\newtheorem{lemma}[equation]{Lemma}
\newtheorem{corollary}[equation]{Corollary}
\newtheorem{conjecture}[equation]{Conjecture}

\theoremstyle{definition}
\newtheorem{example}[equation]{Example}
\newtheorem{definition}[equation]{Definition}

\theoremstyle{remark}
\newtheorem{remark}[equation]{Remark}

\makeatletter\@addtoreset{equation}{section} \makeatother

\author{Ivan Cheltsov}

\title{Log canonical thresholds of del Pezzo surfaces}

\dedicatory{dedicated to Yuri Manin on his seventieth birthday}%

\address{\begin{tabbing}
School of Mathematics\\
University of Edinburgh\\
Edinburgh EH9 3JZ, UK\\
\\
\texttt{I.Cheltsov@ed.ac.uk}
\end{tabbing}}

\begin{document}

\begin{abstract} We study global log canonical thresholds of del Pezzo
surfaces.
\end{abstract}

\maketitle

All varieties are assumed to be defined over $\mathbb{C}$.

\section{Introduction.}
\label{section:intro}

The multiplicity of a nonzero polynomial
$\phi\in\mathbb{C}[z_1,\cdots, z_n]$ at the origin $O\in
\mathbb{C}^n$ is the non\-ne\-ga\-tive integer $m$ such that
$\phi\in\mathfrak{m}^m\setminus\mathfrak{m}^{m+1}$, where
$\mathfrak{m}$ is the maximal ideal of polynomials vanishing at
the point $O$ in $\mathbb{C}[z_1,\cdots, z_n]$. It can be defined
by derivatives, because the equality
$$
m=\mathrm{min}\left\{m\in\mathbb{N}\cup\big\{0\big\}\ \Big|\ \frac{\partial^m \phi\big(z_{1},\ldots,z_{n}\big)}{\partial^{m_1} z_1\partial^{m_2} z_2\cdots\partial^{m_n} z_n}\big(O\big)\ne 0 \right\}.%
$$
holds. We have a similar invariant that is defined by
integrations. This invariant is given by
$$
c_{0}\big(\phi\big)=\mathrm{sup}\left\{c\in\mathbb{Q}\ \Big|\ \text{the function}\ 1\big\slash\big|\phi\big|^{c}~\text{is locally}~L^2~\text{near the point $O\in\mathbb{C}^n$}\right\},%
$$
and $c_{0}(\phi)$~is~called the~log canonical threshold of $\phi$
at the point $O$. The number $c_{0}(\phi)$ appears~in many
places\footnote{The number $c_{0}(\phi)$ is also called the
complex singularity exponent of $\phi$ (see \cite{Ko97}).}. For
instance, it is known that $c_{0}(\phi)$ is the same as the
absolute value of the~largest root of the Bernstein--Sato
polynomial~of~$\phi$ (see \cite{Ko97}).

Even though the log canonical threshold was known implicitly, it
was formally introduced in the paper \cite{Sho93} as follows. Let
$X$ be a variety with log terminal singularities, let $Z\subseteq
X$~be~a~closed subvariety, and let $D$ be an effective
$\mathbb{Q}$-Cartier $\mathbb{Q}$-divisor on $X$. Then the number
$$
\mathrm{lct}_{Z}\big(X,D\big)=\mathrm{sup}\left\{\lambda\in\mathbb{Q}\ \Big|\ \text{the log pair}\ \big(X, \lambda D\big)\ \text{is log canonical  along}~Z\right\}%
$$
is said to be the log canonical threshold of $D$ along $Z$. The
number $\mathrm{lct}_{Z}(X,D)$ is known~to~be~po\-si\-tive and
rational. Moreover, if $X=\mathbb{C}^n$ and $D=(\phi=0)$, then the
equality
$$
\mathrm{lct}_{O}\big(X,D\big)=c_0\big(\phi\big)
$$
holds (see~\cite{Ko97}). For the case $Z=X$ we use the notation
$\mathrm{lct}(X,D)$ instead of $\mathrm{lct}_X(X,D)$. Then
$$
\mathrm{lct}\big(X,D\big)=\mathrm{inf}\left\{\mathrm{lct}_P\big(X,D\big)\ \Big\vert\ P\in X\right\}=\mathrm{sup}\left\{\lambda\in\mathbb{Q}\ \Big|\ \text{the log pair}\ \big(X, \lambda D\big)\ \text{is log canonical}\right\}.%
$$

Even though several methods have been invented in order to compute
log canonical thresholds, it is not easy to compute them in
general. However, the log canonical thresholds play a significant
role in the study on birational geometry  showing many interesting
properties~(see \cite{Ko97}, \cite{Pa01}).

Thus far the log canonical threshold has a local character. In
this paper we wish to develop~its global analogue for Fano
varieties. We shall see it is useful to consider the smallest of
log~cano\-ni\-cal thresholds of effective $\mathbb{Q}$-divisors
numerically equivalent to an anticanonical divisor.

\begin{definition}
\label{definition:threshold} Let $X$ be a Fano variety with  log
terminal singularities, and $G$ be a finite subgroup in
$\mathrm{Aut}(X)$. We define the global $G$-invariant log
canonical threshold of $X$ by the~number
$$
\mathrm{lct}\big(X,G\big)=\mathrm{inf}\left\{\mathrm{lct}(X,D)\ \Big\vert\ \text{the effective $\mathbb{Q}$-divisor}\ D\ \text{is $G$-invariant and}\ D\equiv -K_{X}\right\}.%
$$
\end{definition}

We put $\mathrm{lct}(X)=\mathrm{lct}(X,G)$ if the group $G$ is
trivial. Then it follows from
Definition~\ref{definition:threshold} that
$$
\mathrm{lct}\big(X,G\big)=\mathrm{sup}\left\{\lambda\in\mathbb{Q}\ \left|%
\aligned
&\text{the log pair}\ \left(X, \lambda D\right)\ \text{has log canonical singularities}\\
&\text{for every $G$-invariant effective $\mathbb{Q}$-divisor}\ D\equiv -K_{X}\\
\endaligned\right.\right\}\geqslant 0.%
$$

\begin{example}
\label{example:P-1-1-n} It follows from Proposition~16.9 in
\cite{Ko91} that $\mathrm{lct}(\mathbb{P}(1,1,n))=1/(2+n)$ for
$n\in\mathbb{N}$.
\end{example}

For a given Fano variety, it is usually very hard to compute its
global log canonical threshold explicitly (see \cite{Ch07a}). For
instance, the papers \cite{Hw06} and \cite{Hw06b} show that the
global log canonical~thres\-hold of a~rational homogeneous space
of Picard rank $1$ and Fano index $r$ is $1/r$.

\begin{example}
\label{example:Cheltsov-Park} Let $X$ be a smooth hypersurface in
$\mathbb{P}^{n}$ of degree $n\geq 3$. Then
$$
\mathrm{lct}\big(X\big)\geqslant1-1\big\slash n%
$$
due to \cite{Ch01b}. It is clear that the inequality
$\mathrm{lct}(X)=1-1/n$ holds if the hypersurface $X$ contains a
cone of dimension $n-2$. But the paper \cite{Pu04d} shows that
$\mathrm{lct}(X)=1$ if $X$ is general and $n\geqslant 6$.
\end{example}

Global log canonical thresholds of Fano varieties play an
important role in geometry\footnote{Global log canonical
thresholds of Fano varieties are algebraic counterparts of
$\alpha$-invariants introduced in \cite{Ti87}.}.

\begin{example}
\label{example:CPR} Let $X$ be a general quasismooth hypersurface
in $\mathbb{P}(1,a_{1},\ldots,a_{4})$ of
degree~$\sum_{i=1}^{4}a_{i}$ with terminal singularities such that
$-K_{X}^{3}\leqslant 1$. Then $\mathrm{lct}(X)=1$ by \cite{Ch07a},
which implies that
$$
\mathrm{Bir}\Big(\underbrace{X\times\cdots\times X}_{m\ \text{times}}\Big)=\Big<\prod_{i=1}^{m}\mathrm{Bir}\big(X\big),\ \mathrm{Aut}\Big(\underbrace{X\times\cdots\times X}_{m\ \text{times}}\Big)\Big>,%
$$%
the variety $X\times\cdots\times X$ is not rational and not
birational to a conic bundle (see \cite{Ch07a}).
\end{example}

One of the most interesting applications of global log canonical
thresholds of Fano varieties is the following result proved in
\cite{DeKo01} (see also \cite{Na90} and \cite{Ti87}).

\begin{theorem}
\label{theorem:KE} Let $X$ be a Fano variety with quotient
singularities such that the inequality
$$
\mathrm{lct}\big(X\big)>\frac{\mathrm{dim}\big(X\big)}{\mathrm{dim}\big(X\big)+1}
$$
holds. Then $X$ has an orbifold K\"ahler--Einstein metric.
\end{theorem}

 The following conjecture is inspired by
Question~1~in~\cite{Ti90b}.

\begin{conjecture}
\label{conjecture:stabilization} For a given Fano variety $X$ with
log terminal singularities, the equality
$$
\mathrm{lct}\big(X\big)=\mathrm{lct}\big(X,D\big)
$$
holds for some effective $\mathbb{Q}$-divisor $D$ on the variety
$X$ such that $D\equiv -K_{X}$.
\end{conjecture}

The main purpose of this paper is to prove the following result.

\begin{theorem}
\label{theorem:main}  Let $X$ be a smooth del Pezzo surface. Then
$$
\mathrm{lct}\big(X\big)=\left\{%
\aligned
&1\ \mathrm{when}\ K_{X}^{2}=1\ \mathrm{and}\ |-K_{X}|\ \mathrm{has\ no\ cuspidal\ curves},\\%
&5/6\ \mathrm{when}\ K_{X}^{2}=1\ \mathrm{and}\ |-K_{X}|\ \mathrm{has\ a\ cuspidal\ curve},\\%
&5/6\ \mathrm{when}\ K_{X}^{2}=2\ \mathrm{and}\ |-K_{X}|\ \mathrm{has\ no\ tacnodal\ curves},\\%
&3/4\ \mathrm{when}\ K_{X}^{2}=2\ \mathrm{and}\ |-K_{X}|\ \mathrm{has\ a\ tacnodal\ curve},\\%
&3/4\ \mathrm{when}\ X\ \mathrm{is\ a\ cubic\ surface\ in}\ \mathbb{P}^{3}\ \mathrm{without\ Eckardt\ points},\\%
&2/3\ \mathrm{when}\ K_{X}^{2}=4\ \mathrm{or}\ X\ \mathrm{is\ a\ cubic\ surface\ in}\ \mathbb{P}^{3}\ \mathrm{with\ an\ Eckardt\ point},\\%
&1/2\ \mathrm{when}\ X\cong\mathbb{P}^{1}\times\mathbb{P}^{1}\ \mathrm{or}\ K_{X}^{2}\in\big\{5,6\big\},\\%
&1/3\ \mathrm{in\ the\ remaining\ cases}.\\%
\endaligned\right.%
$$
\end{theorem}

Taking the paper \cite{Pa01} and Theorem~\ref{theorem:main} into
consideration, we see that the assertion of
Conjecture~\ref{conjecture:stabilization} holds for smooth del
Pezzo surfaces. Also, in this paper, we prove the
following~result.

\begin{theorem}
\label{theorem:dP1} Let $X$ be a del Pezzo surface with ordinary double points
such that $K_{X}^{2}=1$.~Then
$$
\mathrm{lct}\big(X\big)=\left\{%
\aligned
&1\ \mathrm{when}\ |-K_{X}|\ \mathrm{does\ not\ have\ cuspidal\ curves},\\%
&3/4\ \mathrm{when}\ |-K_{X}|\ \mathrm{has\ a\ cuspidal\ curve}\ C\ \mathrm{such\ that}\ \mathrm{Sing}(C)\subseteq\mathrm{Sing}(X),\\%
&5/6\ \mathrm{in\ the\ remaining\ cases}.\\%
\endaligned\right.%
$$
\end{theorem}

We see that Theorems~\ref{theorem:KE} and \ref{theorem:dP1} imply
the existence of an~orbifold K\"ahler--Einstein metric on every
del Pezzo surface of degree $1$ that has at most ordinary double
points\footnote{The problem of the existence of a
K\"ahler--Einstein~metric on smooth del Pezzo surfaces is solved
in \cite{Ti90}.}.

\newpage

Further we will study global $G$-invariant log canonical
thresholds of some smooth del Pezzo surfaces admitting an action
of a finite group $G$. Let us consider two examples.

\begin{example}
\label{example:Klein} The simple group
$\mathrm{PGL}(2,\mathrm{F}_{7})$ is a group of automorphisms of
the quartic
$$
x^{3}y+y^{3}z+z^{3}x=0\subset\mathbb{P}^{2}\cong\mathrm{Proj}\Big(\mathbb{C}[x,y,z]\Big),
$$
which induces
$\mathrm{PGL}(2,\mathrm{F}_{7})\subset\mathrm{Aut}(\mathbb{P}^{2})$. Then
$\mathrm{lct}(\mathbb{P}^{2}, \mathrm{PGL}(2,\mathrm{F}_{7}))=4/3$ by
Lemma~\ref{lemma:m-over-r}.
\end{example}

\begin{example}
\label{example:Mukai-Mabuchi} Let $X$ be a del Pezzo surface with
ordinary double points that is given by
$$
\sum_{i=0}^{4}x_{i}^{2}=\sum_{i=0}^{4}\lambda_{i}x_{i}^{2}=0\subseteq\mathbb{P}^{4}\cong\mathrm{Proj}\Big(\mathbb{C}[x_{0},\ldots,x_{4}]\Big),
$$
where $\lambda_{1},\ldots,\lambda_{4}\in\mathbb{C}$. Then
$\mathrm{lct}(X,\mathbb{Z}_{2}^{4})=1$ by
Lemma~\ref{lemma:m-over-r}.
\end{example}

There is a crucial difference between the two and
higher-dimensional cases: in the latter~case, we usually assume
that $G$ is trivial. For surfaces, it is not so, and applications
are more special.

\begin{example}
\label{example:icosahedron} Let $X$ be a smooth cubic surface in
$\mathbb{P}^{3}$ that is given by the equation
$$
x^{2}y+xz^{2}+zt^{2}+tx^{2}=0\subset\mathbb{P}^{3}\cong\mathrm{Proj}\Big(\mathbb{C}[x,y,z,t]\Big),
$$
which implies that $\mathrm{Aut}(X)\cong\mathrm{S}_{5}$ (see \cite{DoIs06}).
Then $\mathrm{lct}(X,\mathrm{S}_{5})=\mathrm{lct}(X,\mathrm{A}_{5})=2$ (see
Lemma~\ref{lemma:m-over-r}), there is a classical embedding
$\mathrm{A}_{5}\subset\mathrm{Aut}(\mathbb{P}^{1})$ and there is a natural
embedding
$$
\mathrm{A}_{5}\subset\mathrm{Aut}\big(\mathbb{P}^{2}\big)\cong\mathrm{PGL}\big(3, \mathbb{C}\big)%
$$
such that the induced embeddings $\mathrm{Aut}(\mathbb{P}^{1}\times
X)\supset\mathrm{A}_{5}\times\mathrm{A}_{5}\subset\mathrm{Aut}(\mathbb{P}^{1}\times\mathbb{P}^{2})$
induce the embeddings
$$
\mathrm{A}_{5}\times\mathrm{A}_{5}\cong\Omega\subset\mathrm{Bir}\big(\mathbb{P}^{3}\big)\supset\Gamma\cong\mathrm{A}_{5}\times\mathrm{A}_{5},
$$
respectively. Then $\Omega$ and $\Gamma$ are not conjugated in
$\mathrm{Bir}(\mathbb{P}^{3})$ by Lemma~\ref{lemma:k-K-square} and
Theorem~\ref{theorem:Pukhlikov}.
\end{example}

We would like to thank H.\,Braden, J.-P.\,Demailly, I.\,Dolgachev,
J.\,Koll\'ar, J.\,Park and Yu.\,Pro\-kho\-rov~for useful comments.
We would like to give a special thanks to J.\,Koll\'ar for
pointing~out a~gap in the~old version of
Theorem~\ref{theorem:Pukhlikov}. We would like to thank the
referee for helpful comments.

\section{Basic tools.}
\label{section:tools}

Let $S$ be a surface with canonical singularities, and $D$ be an
effective $\mathbb{Q}$-divisor on it.

\begin{remark}
\label{remark:convexity} Let $B$ be an effective
$\mathbb{Q}$-divisor on $S$ such that $(S,B)$ is log canonical.
Then
$$
\left(S,\ \frac{1}{1-\alpha}\Big(D-\alpha B\Big)\right)
$$
is not log canonical if $(S, D)$ is not log canonical, where
$\alpha\in\mathbb{Q}$ such that $0\leqslant\alpha<1$.
\end{remark}

Let $\mathrm{LCS}(S,D)\subsetneq S$ be a subset such that
$P\in\mathrm{LCS}(S,D)$ if and only if $(S,D)$ is not log~terminal
at the point $P$. The set $\mathrm{LCS}(S,D)$ is called the locus
of log canonical singularities.

\begin{lemma}
\label{lemma:connectedness} Suppose that $-(K_{S}+D)$ is ample.
Then the set $\mathrm{LCS}(S,D)$ is connected.
\end{lemma}

\begin{proof}
See Theorem~17.4 in \cite{Ko91}.
\end{proof}

Let $P$ be a smooth point of the surface $S$. Suppose that $(S,D)$
is not log canonical at $P$.

\begin{remark}
\label{remark:smooth-points} The inequality
$\mathrm{mult}_{P}(D)>1$ holds (see \cite{Ko97}).
\end{remark}

Let $C$ be an irreducible curve on the surface $S$. Put $D=mC+\Omega$, where
$m$ is a non-negative rational number, and $\Omega$ is an effective
$\mathbb{Q}$-divisor such that $C\not\subseteq\mathrm{Supp}(\Omega)$.

\begin{remark}
\label{remark:curves} Suppose that $C\subseteq\mathrm{LCS}(S,D)$.
Then $m\geqslant 1$ (see \cite{Ko97}).
\end{remark}

Suppose that the inequality $m\leqslant 1$ holds and $P\in C$.

\begin{lemma}
\label{lemma:adjunction} Suppose that $C$ is smooth at $P$. Then
$C\cdot\Omega>1$.
\end{lemma}

\begin{proof}
See Theorem~17.6 in \cite{Ko91}.
\end{proof}

Let $\pi\colon\bar{S}\to S$ be a birational morphism, and $\bar{D}$ is a proper
transform of $D$ via $\pi$. Then
$$
K_{\bar{S}}+\bar{D}+\sum_{i=1}^{r}a_{i}E_{i}\equiv\pi^{*}\big(K_{S}+D\big),
$$
where $E_{i}$ is a $\pi$-exceptional curve, and $a_{i}$ is a
rational number.

\begin{remark}
\label{remark:log-pull-back} The log pair $(S,D)$ is log canonical
if and only if $(\bar{S},\bar{D}+\sum_{i=1}^{r}a_{i}E_{i})$ is log
canonical.
\end{remark}

Suppose that $\pi$ is a blow up of the point $P$. Then $r=1$ and
$\pi(E_{1})=P$. The log pair
$$
\Big(\bar{S},\ \bar{D}+\big(\mathrm{mult}_{P}\big(D\big)-1\big)E_{1}\Big)%
$$
is not log canonical at some point $\bar{P}\in E_{1}$ by
Remark~\ref{remark:log-pull-back}. But
$a_{1}=\mathrm{mult}_{P}(D)-1>0$.

\begin{corollary}
\label{corollary:blow-up-inequality} The inequality
$\mathrm{mult}_{\bar{P}}(\bar{D})+\mathrm{mult}_{P}(D)>2$ holds.
\end{corollary}

Most of the described results are valid in much more general
settings (see \cite{Ko91} and \cite{Ko97}).

\section{Smooth surfaces.}
\label{section:smooth}

In this section we prove Theorem~\ref{theorem:main}. Let $X$ be a
smooth del Pezzo surface. Putting
$$
\omega=\left\{%
\aligned
&1/3\ \mathrm{when}\ X\cong\mathbb{F}_{1}\ \mathrm{or}\ K_{X}^{2}\in\big\{7,9\big\},\\%
&1/2\ \mathrm{when}\ X\cong\mathbb{P}^{1}\times\mathbb{P}^{1}\ \mathrm{or}\ K_{X}^{2}\in\big\{5,6\big\},\\%
&2/3\ \mathrm{when}\ K_{X}^{2}=4\ \mathrm{or}\ X\ \mathrm{is\ a\ cubic\ surface\ in}\ \mathbb{P}^{3}\ \mathrm{with\ an\ Eckardt\ point},\\%
&3/4\ \mathrm{when}\ X\ \mathrm{is\ a\ cubic\ surface\ in}\ \mathbb{P}^{3}\ \mathrm{without\ Eckardt\ points},\\%
&3/4\ \mathrm{when}\ K_{X}^{2}=2\ \mathrm{and}\ |-K_{X}|\ \mathrm{has\ a\ tacnodal\ curve},\\%
&5/6\ \mathrm{when}\ K_{X}^{2}=2\ \mathrm{and}\ |-K_{X}|\ \mathrm{has\ no\ tacnodal\ curves},\\%
&5/6\ \mathrm{when}\ K_{X}^{2}=1\ \mathrm{and}\ |-K_{X}|\ \mathrm{has\ a\ cuspidal\ curve},\\%
&1\ \mathrm{when}\ K_{X}^{2}=1\ \mathrm{and}\ |-K_{X}|\ \mathrm{has\ no\ cuspidal\ curves},\\%
\endaligned\right.%
$$
we see that we must show that $\mathrm{lct}(X)=\omega$ to prove
Theorem~\ref{theorem:main}. But $\mathrm{lct}(X)\leqslant\omega$
by \cite{Pa01}.

Suppose that the inequality $\mathrm{lct}(X)<\omega$ holds. To
prove Theorem~\ref{theorem:main}, we must show that this
assumption leads~to~a~contradiction. There is an effective
$\mathbb{Q}$-divisor $D$ on the surface $X$~such~that the
equivalence $D\equiv -K_{X}$ holds, and $(X, \omega D)$ is not log
canonical at some point $P\in X$.

\begin{lemma}
\label{lemma:degree-1} The inequality $K_{X}^{2}\ne 1$ holds.
\end{lemma}

\begin{proof}
Suppose that $K_{X}^{2}=1$. Take $C\in |-K_{X}|$ such that $P\in
C$. Then $C$ is an~irreducible~curve, and $(X,\omega C)$ is log
canonical. We may assume that $C\not\subseteq\mathrm{Supp}(D)$ by
Remark~\ref{remark:convexity}.~Then
$$
1=C\cdot D\geqslant\mathrm{mult}_{P}\big(D\big)>1/\omega\geqslant 1,%
$$
which is a contradiction. The obtained contradiction completes the
proof.
\end{proof}

\begin{lemma}
\label{lemma:degree-7} The inequality $K_{X}^{2}\leqslant 7$
holds.
\end{lemma}

\begin{proof}
The equalities $\mathrm{lct}(\mathbb{P}^{2})=1/3$ and
$\mathrm{lct}(\mathbb{P}^{1}\times\mathbb{P}^{1})=1/2$ follow from
Remarks~\ref{remark:convexity} and \ref{remark:smooth-points},
which implies that we may assume that $X=\mathbb{F}_{1}$ to
complete the proof. Then $\omega=1/3$.

Let $L$ and $C$ be irreducible curves on $X$ such that $L^{2}=0$
and $C^{2}=-1$. Then
$$
-K_{X}\equiv 2C+3L,
$$
and the singularities of the log pair $(X,\omega(2C+3L))$ are log
canonical.

It follows from Remark~\ref{remark:smooth-points} that
$L\subseteq\mathrm{Supp}(D)$, because $L\cdot D=2$. Therefore, we
may assume that $C\not\subseteq\mathrm{Supp}(D)$ by
Remark~\ref{remark:convexity}. Let $Z$ be a general curve in
$|C+L|$ such that $P\in Z$. Then
$$
3=Z\cdot D\geqslant
\mathrm{mult}_{P}\big(D\big)>1\big\slash\omega=3,
$$
which is a contradiction. The obtained contradiction completes the
proof.
\end{proof}

\begin{lemma}
\label{lemma:degree-4} The inequality $K_{X}^{2}\leqslant 4$
holds.
\end{lemma}

\begin{proof}
Suppose that $K_{X}^{2}\geqslant 5$. Then there is a birational
morphism $\pi\colon X\to S$ such that
\begin{itemize}
\item the morphism $\pi$ is an isomorphism in a~neighborhood of $P$,%
\item either $S\cong\mathbb{F}_{1}$ or $S\cong\mathbb{P}^{1}\times\mathbb{P}^{1}$ or $S\cong\mathbb{P}^{2}$,%
\end{itemize}
and we may assume that $S\cong\mathbb{P}^{1}\times\mathbb{P}^{1}$
whenever $K_{X}^{2}\leqslant 6$. Then the log pair $(S,
\omega\pi(D))$ is not log canonical~at~$\pi(P)$. But $\pi(D)\equiv
-K_{S}$, which is impossible by Lemma~\ref{lemma:degree-7}.
\end{proof}

\begin{lemma}
\label{lemma:degree-4-not} The inequality $K_{X}^{2}\ne 4$ holds.
\end{lemma}

\begin{proof}
Suppose that $K_{X}^{2}=4$. Then $X$ is an intersection of two
quadrics in $\mathbb{P}^{4}$, and
$$
D=\sum_{i=1}^{r}a_{i}C_{i}\equiv -K_{X},
$$
where $C_{i}$ is an irreducible
curve on the surface~$X$, and $0\leqslant a_{i}\in\mathbb{Q}$.

The equality $\omega=2/3$ holds. Suppose that
$a_{k}>1/\omega=3/2$. Then
$$
4=-K_{X}\cdot D=\sum_{i=1}^{r}a_{i}\mathrm{deg}\big(C_{i}\big)\geqslant a_{k}\mathrm{deg}\big(C_{k}\big)>\frac{3\mathrm{deg}\big(C_{k}\big)}{2},%
$$
which implies that $\mathrm{deg}(C_{k})\leqslant 2$. Let $Z$ be an
irreducible curve on $X$ such that $C_{k}+Z$ is cut~out by a
general hyperplane section of $X\subset\mathbb{P}^{4}$ that passes
through $C_{k}$. Then
$$
3\geqslant 4-\mathrm{deg}\big(C_{k}\big)=Z\cdot D=\sum_{i=1}^{r}a_{i}\big(Z\cdot C_{i}\big)\geqslant a_{k}\big(Z\cdot C_{k}\big)=2a_{k}>3,%
$$
which is a contradiction. Therefore, we see that $\omega
a_{i}\leqslant 1$ for every $i=1,\ldots,r$.

There is $\lambda\in\mathbb{Q}$ such that $0<\lambda<\omega=2/3$
and $(X,\lambda D)$ is not log canonical at $P$. Then
$$
\mathrm{LCS}\big(X, \lambda D\big)=\big\{P\big\}
$$
by Lemma~\ref{lemma:connectedness}. But there is a birational
morphism $\pi\colon X\to\mathbb{P}^{2}$ such that $\pi$ is an
isomorphism in a neighborhood of the point $P$. Then $\pi(D)\equiv
-\lambda K_{\mathbb{P}^{2}}$. Let $L$ be a general line on
$\mathbb{P}^{2}$. Then
$$
\pi\big(P\big)\cup L\subseteq\mathrm{LCS}\Big(\mathbb{P}^{2},\ \pi\big(D\big)+L\Big),%
$$
which is impossible by Lemma~\ref{lemma:connectedness}. The
obtained contradiction completes the proof.
\end{proof}

Let $\pi\colon U\to X$ be a blow up of the point $P$, and $E$ be
the exceptional curve of $\pi$. Then
$$
\bar{D}\equiv\pi^{*}\big(-K_{X}\big)-\mathrm{mult}_{P}\big(D\big)E,
$$
where $\bar{D}$ is the proper transform of $D$ on the surface $U$.
It follows from Remark~\ref{remark:log-pull-back} that
$$
\Big(U,\
\omega\bar{D}+\omega\big(\mathrm{mult}_{P}\big(D\big)-1\big)E\Big)
$$
is not log canonical at some point $Q\in E$. Then
$\mathrm{mult}_{Q}(\bar{D})+\mathrm{mult}_{P}(D)>2/\omega$ by
Corollary~\ref{corollary:blow-up-inequality}.

\begin{lemma}
\label{lemma:degree-2} The inequality $K_{X}^{2}\ne 2$ holds.
\end{lemma}

\begin{proof}
Suppose that $K_{X}^{2}=2$. There is a double cover $\psi\colon
X\to\mathbb{P}^{2}$ such that $\psi$ is branched over a smooth
quartic curve $C\subset\mathbb{P}^{2}$. Then either $\psi(P)\in C$
or $\psi(P)\not\in C$.

Suppose that $\psi(P)\in C$. There is a curve $L\in |-K_{X}|$ that
is singular at $P$, and we may~assume that
$L\not\subseteq\mathrm{Supp}(D)$ by Remark~\ref{remark:convexity},
because $(X,\omega L)$ is log canonical. Then
$$
2=L\cdot D\geqslant\mathrm{mult}_{P}\big(D\big)\mathrm{mult}_{P}\big(L\big)\geqslant 2/\omega>2%
$$
in the case when $L$ is irreducible. So, we must have
$L=L_{1}+L_{2}$, where $L_{1}$ and $L_{2}$ are irreducible smooth
curves such that $L_{1}\cdot L_{2}=2$ and
$L_{1}^{2}=L_{2}^{2}=-1$. Put $D=mL_{1}+\Omega$, where $0\leqslant
m\in\mathbb{Q}$, and $\Omega$ is an effective $\mathbb{Q}$-divisor
such that $L_{1}\not\subseteq\mathrm{Supp}(\Omega)$. Then
$$
m+1<2m+\Omega\cdot L_{2}=D\cdot L_{2}=1,%
$$
which is a contradiction. Therefore, we see that $\psi(P)\not\in C$.

In particular, the log pair $(X, \omega D)$ is log canonical outside of
finitely many points.

There is a unique curve $Z\in |-K_{X}|$ such that $P\in Z$ and
$Q\in\bar{Z}$, where $\bar{Z}$ is the proper~transform of the
curve $Z$ on the surface $U$. Then $Z$ consists of at most two
components.

Suppose that $Z$ is irreducible. We may assume
$Z\not\subseteq\mathrm{Supp}(D)$. Hence, we have
$$
2-\mathrm{mult}_{P}\big(D\big)=\bar{Z}\cdot\bar{D}\geqslant\mathrm{mult}_{Q}\big(\bar{D}\big)>2/\omega-\mathrm{mult}_{P}\big(D\big)%
$$
which is a contradiction. So, we must have $Z=Z_{1}+Z_{2}$, where
$Z_{1}$ and $Z_{1}$ are irreducible smooth curves such that
$Z_{1}\cdot Z_{2}=2$ and $Z_{1}^{2}=Z_{2}^{2}=-1$. We may assume
that $P\in Z_{1}$ and $P\not\in Z_{2}$.

It is easy to see that the log pair $(X,\omega Z_{1}+\omega Z_{2})$ is log
canonical. Thus, we may~assume~that either $Z_{1}\not\subseteq\mathrm{Supp}(D)$
or $Z_{2}\not\subseteq\mathrm{Supp}(D)$ by Remark~\ref{remark:convexity}. But
$$
1=Z_{1}\cdot D\geqslant\mathrm{mult}_{P}\big(D\big)\geqslant 1/\omega>1,
$$
which implies that $Z_{2}\not\subseteq\mathrm{Supp}(D)$. Then
$Z_{1}\subseteq\mathrm{Supp}(D)$. Put $D=\bar{m}Z_{1}+\Upsilon$,
where $0<\bar{m}\in\mathbb{Q}$, and $\Upsilon$ is an effective
$\mathbb{Q}$-divisor on the surface $X$ such that
$Z_{1}\not\subseteq\mathrm{Supp}(\Upsilon)$. Then
$$
2\bar{m}\leqslant 2\bar{m}+\Upsilon\cdot Z_{2}=D\cdot Z_{2}=1,%
$$
which gives $\bar{m}\leqslant 1/2$. But $Q\in\bar{Z}_{1}$, where $\bar{Z}$ it
the proper transform of $Z_{1}$ on the surface $U$.~Then
$$
2-\mathrm{mult}_{P}\big(D\big)\geqslant 1-\mathrm{mult}_{P}\big(D\big)+2\bar{m}=\bar{Z}_{1}\cdot\bar{\Upsilon}>2/\omega-\mathrm{mult}_{P}\big(D\big)>2-\mathrm{mult}_{P}\big(D\big)%
$$
by Lemma~\ref{lemma:adjunction}. The obtained contradiction
completes the proof.
\end{proof}

It follows from Lemmas~\ref{lemma:degree-7}, \ref{lemma:degree-4},
\ref{lemma:degree-4-not}, \ref{lemma:degree-1},
\ref{lemma:degree-2} that $X$ is a smooth cubic surface
in~$\mathbb{P}^{3}$.

\begin{lemma}
\label{lemma:cubic-with-Eckardt} The cubic surface $X$ does not
have Eckardt points.
\end{lemma}

\begin{proof}
There is a birational morphism $\pi\colon X\to S$ such that
\begin{itemize}
\item the morphism $\pi$ is an isomorphism in a~neighborhood of the point $P$,%
\item the surface $S$ is a smooth del Pezzo surface and $K_{S}^{2}=4$.%
\end{itemize}

Suppose that $X$ has an Eckardt point\footnote{A point of a cubic
surface is an Eckardt point if the cubic contains $3$ lines
passing through this point.}. Then $\pi(D)\equiv-K_{S}$ and
$(S,\omega\pi(D))$ is not log canonical at the point $\pi(P)$,
which is impossible by Lemma~\ref{lemma:degree-4-not}.
\end{proof}

Therefore, we see that $\omega=3/4$ and $\mathrm{mult}_{P}(D)>4/3$
by Remark~\ref{remark:smooth-points}.

\begin{lemma}
\label{lemma:dp3-point} The log pair $(X,\omega D)$ is log canonical on
$X\setminus P$.
\end{lemma}

\begin{proof}
Arguing as in the proof of Lemma~\ref{lemma:degree-4-not}, we see
that the locus $\mathrm{LCS}(X,\omega D)$ contains~finitely many
points. Then the log pair $(X,\omega D)$ is even log terminal on
$X\setminus P$ by Lemma~\ref{lemma:connectedness}.
\end{proof}

Let $T$ be the unique hyperplane section of $X$ that is singular
at $P$. We may assume that
$$
T\not\subseteq\mathrm{Supp}\big(D\big)
$$
by Remark~\ref{remark:convexity}, because $(S,\omega T)$ is log
canonical. The following cases are possible:
\begin{itemize}
\item the curve $T$ is irreducible and $U$ is a del Pezzo surface;%
\item the curve $T$ is a union of a line and an irreducible conic intersecting at $P$;%
\item the curve $T$ consists of $3$ lines such that one of them does not pass through $P$;%
\end{itemize}
where $T$ is reduced and $-K_{U}$ is nef and big. We exclude these
cases one by one.

\begin{lemma}
\label{lemma:dp3-no-lines} The curve $T$ is reducible.
\end{lemma}

\begin{proof}
Suppose that $T$ is irreducible. There is a double cover $\psi\colon
U\to\mathbb{P}^{2}$ branched over~a~quartic curve. Let $\tau\in\mathrm{Aut}(U)$
be an involution\footnote{The involution $\tau$ induces an involution in
$\mathrm{Bir}(X)$ that is called the Geiser involution.} induced by~$\psi$.~It
follows from \cite{Ma67} that $\tau(\bar{T})=E$~and
$$
\tau^{*}\Big(\pi^{*}\big(-K_{X}\big)\Big)\equiv \pi^{*}\big(-2K_{X}\big)-3E.%
$$

Let $\bar{T}$ be the proper transform of $T$ on the surface $U$. Suppose that
$Q\in\bar{T}$. Then
$$
3-2\mathrm{mult}_{P}\big(D\big)=\bar{T}\cdot\bar{D}\geqslant\mathrm{mult}_{Q}\big(\bar{T}\big)\mathrm{mult}_{Q}\big(\bar{D}\big)>\mathrm{mult}_{Q}\big(\bar{T}\big)\Big(8/3-\mathrm{mult}_{P}\big(D\big)\Big)\geqslant 8/3-\mathrm{mult}_{P}\big(D\big),%
$$
which implies that $\mathrm{mult}_{P}(D)\leqslant 1/3$. But
$\mathrm{mult}_{P}(D)>4/3$. Thus, we see that $Q\not\in\bar{T}$.

Put $\breve{Q}=\pi\circ\tau(Q)$. Let $H$ be the hyperplane section of $X$ that
is singular at $\breve{Q}$. Then $T\ne H$, because $P\ne \breve{Q}$ and $T$ is
smooth outside of the point $P$. Then $P\not\in H$, because otherwise
$$
3=H\cdot T\geqslant\mathrm{mult}_{P}\big(H\big)\mathrm{mult}_{P}\big(T\big)+\mathrm{mult}_{\breve{Q}}\big(H\big)\mathrm{mult}_{\breve{Q}}\big(T\big)\geqslant 4.%
$$

Let $\bar{H}$ be the proper transform of $H$ on the surface $U$. Put
$\bar{R}=\tau(\bar{H})$ and $R=\pi(\bar{R})$. Then
$$
\bar{R}\equiv \pi^{*}\big(-2K_{X}\big)-3E,%
$$
ant the curve $\bar{R}$ must be singular at the point $Q$.

Suppose that $R$ is irreducible. The singularities of the log pair $(X,
\frac{3}{8} R)$ are log canonical, which implies that we may assume that
$R\not\subseteq\mathrm{Supp}(D)$ by Remark~\ref{remark:convexity}. Then
$$
6-3\mathrm{mult}_{P}\big(D\big)=\bar{R}\cdot\bar{D}\geqslant \mathrm{mult}_{Q}(\bar{R})\mathrm{mult}_{Q}\big(\bar{D}\big)>2\Big(8/3-\mathrm{mult}_{P}\big(D\big)\Big),%
$$
which implies that $\mathrm{mult}_{P}(D)<2/3$. But $\mathrm{mult}_{P}(D)>4/3$.
The curve  $R$ must be reducible

The curves $R$ and $H$ are reducible. So, there is a line $L\subset X$ such
that $P\not\in L\ni\breve{Q}$.

Let $\bar{L}$ be the proper transform of $L$ on the surface $U$. Put
$\bar{Z}=\tau(\bar{L})$. Then $\bar{L}\cdot E=0$ and
$$
\bar{L}\cdot\bar{T}=\bar{L}\cdot\pi^{*}(-K_{X})=1,%
$$
which implies that $\bar{Z}\cdot E=1$ and
$\bar{Z}\cdot\pi^{*}\big(-K_{X}\big)=2$. We have $Q\in\bar{Z}$. Then
$$
2-\mathrm{mult}_{P}\big(D\big)=\bar{Z}\cdot\bar{D}\geqslant\mathrm{mult}_{Q}\big(\bar{D}\big)>8/3-\mathrm{mult}_{P}\big(D\big)>2-\mathrm{mult}_{P}\big(D\big)%
$$
in the case when $\bar{Z}\not\subseteq\mathrm{Supp}(\bar{D})$. Hence, we see
that $\bar{Z}\subseteq\mathrm{Supp}(\bar{D})$.

Put $Z=\pi(\bar{Z})$. Then $Z$ is a conic and $P\in Z$. Let $F$ be a line on
$X$ such that $F+Z$~is~cut~out by a hyperplane passing through $Z$. Then
$P\not\in F$, because $T\ne F+Z$.

Put $D=\epsilon Z+\Upsilon$, where $\epsilon$ is a positive rational number,
and $\Upsilon$ is an effective $\mathbb{Q}$-divisor whose support does not
contain the conic $Z$. We may assume that
$F\not\subseteq\mathrm{Supp}(\Upsilon)$ by Remark~\ref{remark:convexity}.~Then
$$
1=F\cdot D=2\epsilon+F\cdot\Upsilon\geqslant 2\epsilon,
$$
which implies that $\epsilon\leqslant 1/2$. Let $\bar{\Upsilon}$ be the proper
transform of $\Upsilon$ on the surface $U$. Then
$$
2-\mathrm{mult}_{P}\big(D\big)+\epsilon=\bar{Z}\cdot\bar{\Upsilon}>8/3-\mathrm{mult}_{P}\big(D\big)%
$$
by Lemma~\ref{lemma:adjunction}, which implies that
$\epsilon>2/3$. But $\epsilon\leqslant 1/2$.
\end{proof}

Therefore, there is a line $L_{1}\subset X$ such that $P\in
L_{1}$.

\begin{lemma}
\label{lemma:dp3-single-line} There is a line $L_{2}\subset X$
such that $L_{1}\ne L_{2}$ and $P\in L_{2}$.
\end{lemma}

\begin{proof} Suppose that there is no line $L_{2}\subset X$ such that $L_{1}\ne L_{2}$ and
$P\in L_{2}$. Then $T=L_{1}+C$, where $C$ is an irreducible conic that passes
through the point $P$.

Let $\bar{L}_{1}$ and $\bar{C}$ be the proper transforms of $L_{1}$ and $C$ on
the surface $U$, respectively. Then
$$
\bar{L}_{1}^{2}=-2,\ -K_{U}\cdot\bar{L}_{1}=0,\ \bar{C}^{2}=-1,\ -K_{U}\cdot\bar{C}=1, %
$$
but the divisor $-K_{U}$ is nef and big. There is a commutative
diagram
$$
\xymatrix{
&U\ar@{->}[d]_{\pi}\ar@{->}[rr]^{\zeta}&&W\ar@{->}[d]^{\psi}\\%
&X\ar@{-->}[rr]_{\rho}&&\mathbb{P}^{2},&}
$$ %
where $\zeta$ is the contraction of the curve $\bar{L}_{1}$ to an
ordinary double point, $\psi$ is a double cover branched over a
quartic curve, and $\rho$ is the projection from the point $P$.

Let $\tau$ be the biregular involution of $U$ induced by~$\psi$.~Then
$\tau(E)=\bar{C}$ and
$$
\tau^{*}\big(\bar{L}_{1}\big)\equiv \bar{L}_{1},\ \tau^{*}\big(E\big)\equiv \bar{C},\ \tau^{*}\Big(\pi^{*}\big(-K_{X}\big)\Big)\equiv \pi^{*}\big(-2K_{X}\big)-3E-\bar{L}_{1}.%
$$

We assume that $T\not\subseteq\mathrm{Supp}(D)$. Then either
$L_{1}\not\subseteq\mathrm{Supp}(D)$ or $C\not\subseteq\mathrm{Supp}(D)$. But
$$
\bar{L}_{1}\cdot\bar{D}=1-\mathrm{mult}_{P}\big(D\big)<0,
$$
which implies that $C\not\subseteq\mathrm{Supp}(D)\supseteq L_{1}$. Put
$D=mL_{1}+\Omega$, where $m$ is a positive rational~number, and $\Omega$ is an
effective $\mathbb{Q}$-divisor whose support does not contain the line $L_{1}$.
Then
$$
m\bar{L}_{1}+\bar{\Omega}\equiv \pi^{*}\big(-K_{X}\big)-\Big(m+\mathrm{mult}_{P}\big(\Omega\big)\Big)E\equiv\pi^{*}\big(-K_{X}\big)-\mathrm{mult}_{P}\big(D\big)E,%
$$
where $\bar{\Omega}$ is the proper transform of $\Omega$ on the surface $U$. We
have
$$
0\leqslant\bar{C}\cdot\bar{\Omega}=2-\mathrm{mult}_{P}(\Omega)+2m<2/3-m,
$$
which implies that $m<2/3$. Then
$\mathrm{mult}_{P}(D)=\mathrm{mult}_{P}(\Omega)+m$, which implies
that
\begin{equation}
\label{equation:single-line}
\mathrm{mult}_{Q}\big(\bar{\Omega}\big)>8/3-\mathrm{mult}_{P}\big(\Omega\big)-m\Big(1+\mathrm{mult}_{Q}\big(\bar{L}_{1}\big)\Big).
\end{equation}

Suppose that $Q\in\bar{L}_{1}$. Then it follows from
Lemma~\ref{lemma:adjunction} that
$$
1-\mathrm{mult}_{P}\big(\Omega\big)+m=\bar{L}_{1}\cdot\bar{\Omega}>8/3-\mathrm{mult}_{P}\big(\Omega\big)-m,%
$$
which implies that $m>5/6$. But $m<2/3$. Hence, we see that
$Q\not\in\bar{L}_{1}$.

Suppose that $Q\in\bar{C}$. Then it follows from the
inequality~\ref{equation:single-line} that
$$
2-\mathrm{mult}_{P}\big(\Omega\big)-2m=\bar{C}\cdot\bar{\Omega}>8/3-\mathrm{mult}_{P}\big(\Omega\big)-m,%
$$
which implies that $m<0$. Hence, we see that $Q\not\in\bar{C}$.

We have $\tau(E)=\bar{C}$. Let $H$ be the hyperplane section of
the cubic surface $X$ that is singular at the point
$\pi\circ\tau(Q)\in C$. Then $P\not\in H$, because $C$ is smooth.

Let $\bar{H}$ be the proper transform of $H$ on the surface $U$. Put
$\bar{R}=\tau(\bar{H})$ and $R=\pi(\bar{R})$. Then
$$
\bar{R}\equiv \pi^{*}\big(-2K_{X}\big)-3E-\bar{L}_{1},%
$$
and the curve $\bar{R}$ is singular at the point $Q$ by construction.

Suppose that $R$ is irreducible. Then $R+L_{1}\equiv -2K_{X}$, but
$(X,\frac{3}{8}(R+L_{1}))$ is log canonical, which implies that we may assume
that $R\not\subseteq\mathrm{Supp}(D)$ by Remark~\ref{remark:convexity}. The
inequality~\ref{equation:single-line}~gives
$$
5-2\Big(m+\mathrm{mult}_{P}\big(\Omega\big)\Big)-m=\bar{R}\cdot\bar{\Omega}\geqslant 2\mathrm{mult}_{Q}\big(\bar{\Omega}\big)>2\Big(8/3-m-\mathrm{mult}_{P}\big(\Omega\big)\Big),%
$$
which implies that $m<0$. Hence, there is a line $L\subset X$ such
that $P\not\in L$ and $\pi\circ\tau(Q)\in L$.

Let $\bar{L}$ be the proper transform of the line $L$ on the
surface $U$. Then
$$
\bar{L}\cdot\bar{C}=\bar{L}\cdot\pi^{*}(-K_{X})=1\ \mathrm{and}\ \bar{L}\cdot E=\bar{L}\cdot\bar{L}_{1}=0,%
$$
but $\tau$ preserves the intersection form. Put
$\bar{Z}=\tau(\bar{L})$. Then $\bar{Z}\cdot E=1$,
$\bar{Z}\cdot\bar{L}_{1}=0$, $\bar{Z}\cdot\pi^{*}(-K_{X})=2$.

Suppose that the support of $\bar{\Omega}$ does not contain
$\bar{Z}$. Then the inequality~\ref{equation:single-line}
implies~that
$$
2-m-\mathrm{mult}_{P}\big(\Omega\big)=\bar{Z}\cdot\bar{\Omega}>8/3-m-\mathrm{mult}_{P}\big(\Omega\big),
$$
which is impossible. Thus, the support of $\bar{\Omega}$ must
contain the curve $\bar{Z}$.

Put $Z=\pi(\bar{Z})$. Then $Z$ is a conic that passes through the
point $P$. Let $F$ be a line on $X$~such that the curve $F+Z$ is
cut out by a hyperplane passing through $Z$. Then $P\not\in F$.
Put
$$
D=\epsilon Z+mL_{1}+\Upsilon,
$$
where $\epsilon$ is a positive rational number, and $\Upsilon$ is
an effective $\mathbb{Q}$-divisor on the surface $X$ such that the
support of the divisor $\Upsilon$ does not contain the curves $Z$
and $L_{1}$.

We may assume that $F\not\subseteq\mathrm{Supp}(\Upsilon)$, because $(X, \omega
(F+Z))$ is log canonical.~Then
$$
1=F\cdot D=2\epsilon+mF\cdot L_{1}+F\cdot\Upsilon\geqslant 0,
$$
which implies that $\epsilon\leqslant 1/2$. But
$\bar{Z}\cap\bar{L}_{1}=\varnothing$. Then it follows from
Lemma~\ref{lemma:adjunction} that
$$
2-\mathrm{mult}_{P}\big(D\big)+\epsilon=\bar{Z}\cdot\bar{\Upsilon}>8/3-\mathrm{mult}_{P}\big(D\big),%
$$
where $\bar{\Upsilon}$ is a proper transform of $\Upsilon$ on the surface
$U$. We deduce that  $\epsilon>2/3$. But $\epsilon\leqslant 1/2$.
\end{proof}

We have $T=L_{1}+L_{2}+L_{3}$, where $L_{3}$ is a line such that $P\not\in
L_{3}$. Then
$$
\bar{L}_{1}^{2}=\bar{L}_{2}^{2}=-2,\
E\cdot\bar{L}_{1}=E\cdot\bar{L}_{2}=-K_{U}\cdot\bar{L}_{3}=1,\
-K_{U}\cdot\bar{L}_{1}=-K_{U}\cdot\bar{L}_{2}=E\cdot\bar{L}_{3}=0,\
\bar{L}_{3}^{2}=-1,%
$$
where $\bar{L}_{i}$ is the proper transform of $L_{i}$ on the surface $U$.
There is a commutative diagram
$$
\xymatrix{
&U\ar@{->}[d]_{\pi}\ar@{->}[rr]^{\zeta}&&W\ar@{->}[d]^{\psi}\\%
&X\ar@{-->}[rr]_{\rho}&&\mathbb{P}^{2},&}
$$ %
where $\zeta$ is the contraction of the curves $\bar{L}_{1}$ and
$\bar{L}_{2}$ to ordinary double points, $\psi$ is a double cover
branched over a quartic curve, and $\rho$ is the projection from
the point $P$.

Let $\tau$ be the biregular involution of the surface $U$ induced
by~$\psi$.~Then
$$
\tau^{*}\Big(\pi^{*}\big(-K_{X}\big)\Big)\equiv \pi^{*}\big(-2K_{X}\big)-3E-\bar{L}_{1}-\bar{L}_{2},%
$$
and $\tau(\bar{L}_{1})=\bar{L}_{1}$, $\tau(\bar{L}_{2})=\bar{L}_{2}$,
$\tau(\bar{L}_{3})=E$. Recall that $\mathrm{mult}_{P}(D)>4/3$ by
Remark~\ref{remark:smooth-points}.

We assume that $T\not\subseteq\mathrm{Supp}(D)$. Then $\mathrm{Supp}(D)$ does
not contain one of $L_{1}$, $L_{2}$, $L_{3}$. But
$$
\bar{L}_{1}\cdot\bar{D}=\bar{L}_{2}\cdot\bar{D}=1-\mathrm{mult}_{P}\big(D\big)<0,
$$
which implies that $L_{2}\subseteq\mathrm{Supp}(D)\supseteq L_{2}$
and $L_{3}\not\subseteq\mathrm{Supp}(D)$. Put
$$
D=m_{1}L_{1}+m_{2}L_{2}+\Omega,
$$
where $0<m_{i}\in\mathbb{Q}$, and $\Omega$ is an effective
$\mathbb{Q}$-divisor such that
$L_{2}\not\subseteq\mathrm{Supp}(\Omega)\not\supseteq L_{2}$.

The inequality $m_{1}+m_{2}\leqslant 1$ holds, because
$1-m_{1}-m_{2}=L_{3}\cdot\Omega\geqslant 0$.

Let $\bar{\Omega}$ be the proper transform of $\Omega$ on the surface $U$. Then
$$
m_{1}\bar{L}_{1}+m_{2}\bar{L}_{2}+\bar{\Omega}\equiv \pi^{*}\big(-K_{X}\big)-\Big(m_{1}+m_{2}+\mathrm{mult}_{P}\big(\Omega\big)\Big)E,%
$$
where
$m_{1}+m_{2}+\mathrm{mult}_{P}(\Omega)=\mathrm{mult}_{P}(D)$. The
latter equality implies that
\begin{equation}
\label{equation:two-lines}
\mathrm{mult}_{Q}\big(\bar{\Omega}\big)>8/3-\mathrm{mult}_{P}\big(\Omega\big)-m_{1}\Big(1+\mathrm{mult}_{Q}\big(\bar{L}_{1}\big)\Big)-m_{1}\Big(1+\mathrm{mult}_{Q}\big(\bar{L}_{2}\big)\Big).
\end{equation}

\begin{lemma}
\label{lemma:cubic-upstairs-line} The curves $\bar{L}_{1}$ and
$\bar{L}_{2}$ do not contain the point $Q$.
\end{lemma}

\begin{proof}
Suppose that $Q\in\bar{L}_{1}\cup\bar{L}_{2}$. Without loss of generality we
may assume that $Q\in\bar{L}_{1}$. Then
$$
1-\mathrm{mult}_{P}\big(\Omega\big)-m_{2}+m_{1}=\bar{L}_{1}\cdot\bar{\Omega}>8/3-\mathrm{mult}_{P}\big(\Omega\big)-m_{1}-m_{2}%
$$
by Lemma~\ref{lemma:adjunction}. We have $m_{1}>5/6$. Then
$$
1-m_{1}+m_{2}=\Omega\cdot L_{2}>4/3-m_{1}-m_{2},
$$
which implies the inequality $m_{2}>1/6$. The latter contradicts
the inequality  $m_{1}+m_{2}\leqslant 1$.
\end{proof}

Therefore, the point $\pi\circ\tau(Q)$ is contained in the line
$L_{3}$, but $\pi\circ\tau(Q)\not\in L_{1}\cup L_{2}$ .

\begin{lemma}
\label{lemma:cubic-no-second-line} The line $L_{3}$ is the only
line on $X$ that passes through the point $\pi\circ\tau(Q)$.
\end{lemma}

\begin{proof}
Suppose that there is a line $L\subset X$ such that $L\ne L_{3}$
and $\pi\circ\tau(Q)\in L$.  Then
$$
\bar{L}\cdot\bar{L}_{1}=\bar{L}\cdot\bar{L}_{2}=\bar{L}\cdot E=0,\ \bar{L}\cdot\pi^{*}(-K_{X})=\bar{L}\cdot\bar{L}_{3}=1,%
$$
where $\bar{L}$ is the proper transform of the line $L$ on the
surface $U$.

The involution $\tau$ preserves the intersection form. Put
$\bar{Z}=\tau(\bar{L})$ and $Z=\pi(\bar{Z})$. Then
$$
\bar{Z}\cdot E=1,\ \bar{Z}\cdot\bar{L}_{3}=0,\ \bar{Z}\cdot\pi^{*}\big(-K_{X}\big)=2,%
$$
which implies that the curve $\pi(\bar{Z})$ is a conic passing
through the point $P$.

The support of the divisor $\Omega$ contains the conic $Z$, because otherwise
$$
2-m_{1}-m_{2}-\mathrm{mult}_{P}\big(\Omega\big)=\bar{Z}\cdot\bar{\Omega}>8/3-m_{1}-m_{2}-\mathrm{mult}_{P}\big(\Omega\big),
$$
which is impossible. Put $D=\epsilon Z+m_{1}L_{1}+m_{2}L_{2}+\Upsilon$, where
$\epsilon$ is a positive rational number, and $\Upsilon$ is an effective
$\mathbb{Q}$-divisor on $X$ whose support does not contain
$Z$,~$L_{1}$,~$L_{2}$.

Let $F$ be a line on the surface $X$ such that the curve $F+Z$ is
cut out by a hyperplane that passes through the curve $Z$. We may
assume that the~support of $\Upsilon$ does not contain $F$. Then
$$
1=F\cdot D=2\epsilon+m_{1}F\cdot L_{1}+m_{2}F\cdot L_{2}+F\cdot\Upsilon\geqslant 0,%
$$
which implies that $\epsilon\leqslant 1/2$. But $Q\not\in\bar{L}_{1}$ and
$Q\not\in\bar{L}_{2}$ by Lemma~\ref{lemma:cubic-upstairs-line}. Thus, the log
pair
$$
\Big(U,\ \epsilon\bar{Z}+\omega \bar{\Upsilon}+\big(\omega\mathrm{mult}_{P}\big(D\big)-1\big)E\Big)%
$$
is not log canonical at the point $Q$, where $\bar{\Upsilon}$ is a proper
transform of $\Upsilon$ on the surface $U$. Then
$$
2-\mathrm{mult}_{P}\big(D\big)+\epsilon=2-\mathrm{mult}_{P}\big(D\big)+\epsilon-m_{1}\bar{L}_{1}\cdot\bar{Z}-m_{2}\bar{L}_{2}\cdot\bar{Z}=\bar{Z}\cdot\bar{\Upsilon}>8/3-\mathrm{mult}_{P}\big(D\big)%
$$
by Lemma~\ref{lemma:adjunction}, which implies that
$\epsilon>2/3$. But $\epsilon\leqslant 1/2$.
\end{proof}

Let $C\subset X$ be a  conic such that $C+L_{3}$ is cut out by the
hyperplane tangent to $X$ at $\pi\circ\tau(Q)$, and let $\bar{C}$
be the proper transform of $C$ on the surface $U$. Put
$\bar{Z}=\tau(\bar{C})$ and $Z=\pi(\bar{Z})$. Then
$$
\bar{Z}\equiv \pi^{*}\big(-2K_{X}\big)-4E-\bar{L}_{1}-\bar{L}_{2},
$$
and $Z$ is singular at $P$. We have $\bar{Z}\cdot E=2$ and
$\bar{Z}\cdot\bar{L}_{1}=\bar{Z}\cdot\bar{L}_{2}=0$, because
$C\cap L_{1}=C\cap L_{2}=\varnothing$.

\begin{lemma}
\label{lemma:cubic-final} The support of the divisor $D$ contains
$Z$.
\end{lemma}

\begin{proof}
Suppose that $Z\not\subseteq\mathrm{Supp}(D)$. Then it follows
from Corollary~\ref{corollary:blow-up-inequality} that
$$
4-2\mathrm{mult}_{P}\big(D\big)=\bar{Z}\cdot\bar{D}\geqslant\mathrm{mult}_{Q}\big(\bar{D}\big)>8/3-\mathrm{mult}_{P}\big(D\big),%
$$
which implies that $\mathrm{mult}_{P}(D)<4/3$. But $\mathrm{mult}_{P}(D)>4/3$.
\end{proof}

Put $D=\epsilon Z+m_{1}L_{1}+m_{2}L_{2}+\Upsilon$, where
$0<\epsilon\in\mathbb{Q}$, and $\Upsilon$ is an effective
$\mathbb{Q}$-divisor whose support does not contain the curves
$Z$, $L_{1}$, $L_{2}$. Then $L_{1}+L_{2}+Z\equiv -2K_{X}$ and
$$
D\cdot L_{1}=m_{2}-m_{1}+2\epsilon+L_{1}\cdot\Upsilon=D\cdot L_{2}=m_{1}-m_{2}+2\epsilon+L_{2}\cdot\Upsilon=1,%
$$
which implies that $\epsilon\leqslant 1/2$. Let $\bar{\Upsilon}$ be a proper
transform of $\Upsilon$ on the surface $U$. Then
$$
4-2\mathrm{mult}_{P}\big(D\big)=\bar{Z}\cdot\bar{\Upsilon}>8/3-\mathrm{mult}_{P}\big(D\big)%
$$
by Lemma~\ref{lemma:adjunction}, which implies that
$\mathrm{mult}_{P}(D)<4/3$. But $\mathrm{mult}_{P}(D)>4/3$.

The obtained contradiction completes the proof Theorem~\ref{theorem:main}.

\section{Singular surfaces.}
\label{section:dP-1-2-singular}

Let $X$ be a del Pezzo surface with Du Val singularities such that
$K_{X}^{2}=1$, and singularities of the surface $X$ consist of
finitely many points of type $\mathbb{A}_{1}$ or $\mathbb{A}_{2}$.
Put
$$
\omega=\left\{%
\aligned
&1\ \text{when}\ |-K_{X}|\ \text{does not have cuspidal curves},\\%
&2/3\ \text{when}\ |-K_{X}|\ \text{has a cuspidal curve}\ C\ \text{such that}\ \text{Sing}(C)\ \text{is a point of type}\ \mathbb{A}_{2},\\%
&5/6\ \text{when}\ |-K_{X}|\ \text{has cuspidal curves, but their cusps are not contained in}\ \text{Sing}(S),\\%
&3/4\ \text{in\ the\ remaining\ cases}.\\%
\endaligned\right.%
$$

\begin{lemma}
\label{lemma:singular} The equality $\mathrm{lct}(X)=\omega$
holds.
\end{lemma}

\begin{proof}
Taking into a consideration curves in $|-K_{X}|$, we see that
$\mathrm{lct}(X)\leqslant\omega$. Thus,~to~conclude the proof, we
may assume that $\mathrm{lct}(X)<\omega$. Then there is an
effective $\mathbb{Q}$-divisor $D$ on the~sur\-face~$X$ such~that
$D\equiv-K_{X}$, but $(X, \lambda D)$ is not log terminal and for
some $\omega>\lambda\in\mathbb{Q}$.

Suppose that $\mathrm{LCS}(X,\lambda D)$ is not zero-dimensional.
There is an irreducible curve $C$ such that
$$
D=mC+\Omega
$$
where $1<1/\lambda\leqslant m\in\mathbb{Q}$, and $\Omega$ is an
effective $\mathbb{Q}$-divisor such that
$C\not\subseteq\mathrm{Supp}(\Omega)$. Then
$$
1=H\cdot D=mH\cdot C+H\cdot\Omega>m>1,
$$
where $H$ is a general curve in the pencil $|-K_{X}|$. Thus, the
locus $\mathrm{LCS}(X,\lambda D)$ is zero-dimensional.

It follows from Lemma~\ref{lemma:connectedness} that the locus
$\mathrm{LCS}(X,\lambda D)$ consists of a single point $P\in X$.

Let $Z$ be the curve in $|-K_{X}|$ such that $P\in Z$. Arguing as
in the proof of Lemma~\ref{lemma:degree-1},~we~see that we may
assume that $P\in\mathrm{Sing}(X)$.

We may assume that $Z\not\subseteq\mathrm{Supp}(D)$, because $(X,\omega Z)$ is
log canonical, and $Z$ is irreducible.

Suppose that $P$ is a point of type $\mathbb{A}_{1}$. Let
$\pi\colon U\to X$ be a blow up of the point $P$. Then
$$
\left\{\aligned
&\bar{D}\equiv \pi^{*}\big(-K_{X}\big)-aE,\\%
&\bar{Z}\equiv\pi^{*}\big(-K_{X}\big)-E,\\%
\endaligned
\right.
$$
where $\bar{D}$ and $\bar{Z}$ are proper transforms of  $D$ and
$Z$ on the surface $U$, respectively, $E$ is the~exceptional curve
of $\pi$, and $a$ is a positive rational number. Then $a\leqslant
1/2$, because $1-2a=\bar{Z}\cdot\bar{D}\geqslant 0$.

The log pair$(U, \lambda\bar{D}+\lambda a E)$ is not log terminal
at some point $Q\in E$ by Remark~\ref{remark:log-pull-back}. Then
$$
1\geqslant 2a=E\cdot\bar{D}>1/\lambda>1%
$$
by Lemma~\ref{lemma:adjunction}, which is a contradiction. Thus,
the point $P$ is a singular point of type $\mathbb{A}_{2}$.

There is a birational morphism $\zeta\colon W\to X$ such that
$\zeta$ contracts two irreducible smooth rational curves $E_{1}$
and $E_{2}$ to the point $P$, the morphism $\zeta$ induces an
isomorphism
$$
W\setminus\Big(E_{1}\cup E_{2}\Big)\cong X\setminus P,
$$
and $W$ is smooth along $E_{1}$ and $E_{2}$. Then
$E_{1}^{2}=E_{2}^{2}=-2$ and $E_{1}\cdot E_{2}=1$. But
$$
\left\{\aligned
&\grave{D}\equiv \zeta^{*}\big(-K_{X}\big)-a_{1}E_{1}-a_{2}E_{2},\\%
&\grave{Z}\equiv\zeta^{*}\big(-K_{X}\big)-E_{1}-E_{2}E,\\%
\endaligned
\right.
$$
where $\grave{D}$ and $\grave{Z}$ are proper transforms of $D$ and
$Z$ on the surface $W$, respectively, and $0\leqslant
a_{i}\in\mathbb{Q}$.

The inequalities $\grave{Z}\cdot\grave{D}\geqslant 0$,
$E_{1}\cdot\grave{D}\geqslant 0$, $E_{1}\cdot\grave{D}\geqslant 0$
imply that
$$
a_{1}+a_{2}\leqslant 1,\ 2a_{1}\geqslant a_{2},\ 2a_{2}\geqslant a_{1},%
$$
respectively. Thus, we see that $a_{1}\leqslant 2/3$ and $a_{2}\leqslant 2/3$.
But the equivalence
$$
K_{W}+\lambda\grave{D}+\lambda a_{1}E_{1}+\lambda a_{2}E_{2}\equiv \zeta^{*}\Big(K_{X}+\lambda D\Big)%
$$
implies the existence of a point $O\in E_{1}\cup E_{2}$ such that
$(W,\lambda\grave{D}+\lambda a_{1}E_{1}+\lambda a_{2}E_{2})$ is not
log~ter\-mi\-nal at the point $O$ (see Remark~\ref{remark:log-pull-back}).
Without loss of generality, we may assume that $O\in E_{1}$.

Suppose that $O\not\in E_{2}$. Then $(W, \lambda\grave{D}+E_{1})$
is not log terminal at $Q$. We have
$$
2a_{1}-a_{2}=E_{1}\cdot\grave{D}>1/\lambda>1,%
$$
by Lemma~\ref{lemma:adjunction}, which implies that $a_{1}>2/3$,
because $2a_{2}\geqslant a_{1}$. But $a_{1}\leqslant 2/3$.

Thus, we see that $O=E_{1}\cap E_{2}$. Then
$$
\left\{%
\aligned
&2a_{1}-a_{2}=E_{1}\cdot\grave{D}\geqslant 1/\lambda-a_{2}>1-a_{2},\\%
&2a_{2}-a_{1}=E_{1}\cdot\grave{D}\geqslant 1/\lambda-a_{1}>1-a_{1},\\%
\endaligned\right.%
$$
by Lemma~\ref{lemma:adjunction}, which implies that $a_{1}>1/2$
and $a_{2}>1/2$. But $a_{1}+a_{2}\leqslant 1$.
\end{proof}

The assertion of Theorem~\ref{theorem:dP1} follows from
Lemma~\ref{lemma:singular}.

\section{Invariant thresholds.}
\label{section:invariant}

Let $X$ is a smooth del Pezzo surface, let $H$ be a Cartier divisor on $X$, let
$G$ be a finite subgroup in $\mathrm{Aut}(X)$ such that~the~$G$-in\-va\-ri\-ant
subgroup of the group $\mathrm{Pic}(X)$ is $\mathbb{Z}H$, and
\begin{itemize}
\item let $r$ be the biggest natural number such that $-K_{X}\sim rH$,%
\item let $k$ be the smallest natural number such that $k=|\Sigma|$, where $\Sigma\subset X$ is a $G$-orbit,%
\item let $m$ be the smallest natural number such that there is a $G$-invariant divisor in $|mH|$.%
\end{itemize}

It follows from Definition~\ref{definition:threshold} that $\mathrm{lct}(X,
G)\leqslant m/r$.

\begin{lemma}
\label{lemma:m-over-r} Suppose that
$h^{0}(X,\mathcal{O}_{X}((m-r)H))<k$. Then
$\mathrm{lct}(X,G)=m/r$.
\end{lemma}

\begin{proof}
We suppose that $\mathrm{lct}(X,G)<m/r$. Then there is an
effective $G$-invariant $\mathbb{Q}$-divisor $D$ on the~surface
$X$ such that $\mathrm{LCS}(X,\lambda D)\ne\varnothing$  and
$D\equiv-K_{X}$, where $0<\lambda\in\mathbb{Q}$ such that
$\lambda<m/r$.

It follows from the Nadel vanishing theorem (see Theorem~9.4.8 in
\cite{La04}) that the sequence
\begin{equation}
\label{equation:m-over-r}
H^{0}\Big(X,\ \mathcal{O}_{X}\big((m-r)H\big)\Big)\longrightarrow H^{0}\Big(\mathcal{O}_{\mathcal{L}}\otimes\mathcal{O}_{X}\big((m-r)H\big)\Big)\longrightarrow 0%
\end{equation}
is exact, where $\mathcal{J}(\lambda D)$ is the multiplier ideal
sheaf of $\lambda D$, and $\mathcal{L}$ is the corresponding
subscheme.

Suppose that $\mathcal{L}$ is zero-dimensional. Then the exact
sequence~\ref{equation:m-over-r} implies that
$$
k>h^{0}\Big(X,\ \mathcal{O}_{X}\big((m-r)H\big)\Big)\geqslant h^{0}\Big(\mathcal{O}_{\mathcal{L}}\otimes\mathcal{O}_{X}\big((m-r)H\big)\Big)=h^{0}\big(\mathcal{O}_{\mathcal{L}}\big)\geqslant\Big\vert\mathrm{Supp}\big(\mathcal{L}\big)\Big\vert\geqslant k,%
$$
because the subscheme $\mathcal{L}$ is $G$-invariant. Hence, the
subscheme $\mathcal{L}$ is not zero-dimensional.

Thus, there is a $G$-invariant reduced curve $C$ on the surface
$X$ such that
$$
\lambda D=\mu C+\Omega,
$$
where $\mu\geqslant 1$, and $\Omega$ is an effective one-cycle on
the surface $X$, whose support does not contain any component of
the curve $C$. Then $C\sim l H$ for some natural number $l$. We
have $l\geqslant m$. But
$$
m>\lambda r\geqslant \mu l\geqslant l\geqslant m,%
$$
because the~$G$-in\-va\-ri\-ant subgroup of the group
$\mathrm{Pic}(X)$ is generated by the divisor $H$.
\end{proof}

Let us show how to apply Lemma~\ref{lemma:m-over-r}.

\begin{example}
\label{example:dP6} Suppose that $K_{X}^{2}=6$ and $k\ne 1$. Then
$X$ has $6$ curves $E_{1},\ldots,E_{6}$ such that
$$
\sum_{i=1}^{6}E_{i}\sim -K_{X}
$$
and $E_{i}^{2}=-1$. The divisor $\sum_{i=1}^{6}E_{i}$ is
$G$-invariant. Then $\mathrm{lct}(X, G)=1$ by
Lemma~\ref{lemma:m-over-r}.
\end{example}

\begin{example}
\label{example:dP9} Suppose that $X=\mathbb{P}^{2}$ and $G=\mathrm{A}_{5}$ such
that the subgroup $G$ leaves invariant a~smooth conic on $\mathbb{P}^{2}$. Then
$\mathrm{lct}(X,G)=2/3$ by Lemma~\ref{lemma:m-over-r}, because $r=3$, $k=6$,
$m=2$.
\end{example}

\begin{example}
\label{example:quintic} Suppose that $K_{X}^{2}=6$ and
$G=\mathrm{Aut}(X)\cong\mathrm{S}_{5}$ (see \cite{RaSl02}). Then $r=1$
and~$k>6$, because the stabilizer of every point induces a faithful
two-dimensional linear representation in its~tangent space. Then
$\mathrm{lct}(X,G)=2$ by Lemma~\ref{lemma:m-over-r}, because $m=2$ (see
\cite{RaSl02}).
\end{example}

Even if $h^{0}(X,\mathcal{O}_{X}((m-r)H))\geqslant k$, we still may be able to
show that $\mathrm{lct}(X,G)=m/r$.

\begin{lemma}
\label{lemma:Fermat-cubic} Suppose that $X$ be the cubic surface
in $\mathbb{P}^{3}$ that is given by the equation
$$
x^{3}+y^{3}+z^{3}+t^{3}=0\subset\mathbb{P}^{3}\cong\mathrm{Proj}\Big(\mathbb{C}[x,y,z,t]\Big),
$$
and $G=\mathrm{Aut}(X)$. Then  $\mathrm{lct}(X, G)=4$
\end{lemma}

\begin{proof}
We have $r=1$ and $G\cong\mathbb{Z}_{3}^{3}\rtimes\mathrm{S}_{4}$
(see \cite{DoIs06}). Then it is easy to check that $m=4$ and
$k=18$, which implies that we are unable to apply
Lemma~\ref{lemma:m-over-r} to deduce the equality
$\mathrm{lct}(X,G)=4$.

Suppose that $\mathrm{lct}(X,G)<4$. Then there is an effective
$G$-invariant $\mathbb{Q}$-divisor $D$ on the~cubic surface $X$
such that $\mathrm{LCS}(X,\lambda D)\ne\varnothing$  and
$D\equiv-K_{X}$, where $0<\lambda\in\mathbb{Q}$ such that
$\lambda<4$.

Arguing as in the proof of Lemma~\ref{lemma:m-over-r}, we see that
the locus $\mathrm{LCS}(X,\lambda D)$ consists of $18$ points,
because every $G$-orbit containing at most $20$ points must
consist of $18$ points. Then
$$
\mathrm{LCS}\big(X,\ \lambda D\big)=\big\{O_{1},\ldots,O_{18}\big\},%
$$
where $O_{1},\ldots,O_{18}$ are all Eckardt points of the surface
$X$ (see \cite{DoIs06}).

Let $R$ be a curve on the surface $X$ that is cut out by $xyzt=0$.
Then $R$ is~$G$-invariant, and the log pair $(X, R)$ is log
canonical. We may assume that $R\not\subseteq\mathrm{Supp}(D)$ by
Remark~\ref{remark:convexity}. Then
$$
12=R\cdot
D\geqslant\sum_{i=1}^{18}\mathrm{mult}_{O_{i}}\big(R\big)\mathrm{mult}_{O_{i}}\big(D\big)=\sum_{i=1}^{18}2\mathrm{mult}_{O_{i}}\big(D\big)\geqslant 36\mathrm{mult}_{O_{i}}\big(D\big),%
$$
which implies that $\mathrm{mult}_{O_{i}}(D)\leqslant 1/3$.

Let $\pi\colon U\to X$ be a blow up of the points
$O_{1},\ldots,O_{18}$. Then
$$
K_{U}+4\bar{D}+\sum_{i=1}^{18}\Big(4\mathrm{mult}_{O_{i}}\big(D\big)-1\Big)E_{i}\equiv \pi^{*}\Big(K_{X}+4D\Big),%
$$
where $E_{i}$ is the $\pi$-exceptional curve such that $\pi(E_{i})=O_{i}$, and
$\bar{D}$ is the proper transform of $D$ on the surface $U$. Then there is
$Q_{i}\in E_{i}$ such that
$\mathrm{mult}_{Q_{i}}(\bar{D})>1/2-\mathrm{mult}_{O_{i}}(D)$ for
$i=1,\ldots,18$.

Let $\Sigma$ be the $G$-orbit of the point $Q_{i}$. Then
$\Sigma\cap E_{i}\ne Q_{i}$, because the representation induced by
the action of the stabilizer of $O_{i}$ on its tangent space is
irreducible. We have
$$
\mathrm{mult}_{O_{i}}\big(D\big)=E_{i}\cdot\bar{D}>\Big|\Sigma\cap E_{i}\Big|\Big(1/2-\mathrm{mult}_{O_{i}}\big(D\big)\Big),%
$$
which implies that $|\Sigma\cap E_{i}|=1$, because
$\mathrm{mult}_{O_{i}}(D)\leqslant 1/3$.
\end{proof}

\begin{lemma}
\label{lemma:quintic-A5} Suppose that $K_{X}^{2}=5$ and
$G=\mathrm{A}_{5}$. Then $\mathrm{lct}(X,G)=2$.
\end{lemma}

\begin{proof}
The surface $X$ is embedded in $\mathbb{P}^{5}$ by the linear system
$|-K_{X}|$, and $X$~contains $10$ lines, which we denote as
$L_{1},\ldots,L_{10}$. Then $r=1$ and $\mathrm{Aut}(X)\cong\mathrm{S}_{5}$ (see
\cite{RaSl02}).

The divisor $\sum_{i=1}^{10}L_{i}\sim -2K_{X}$ is $\mathrm{S}_{5}$-invariant,
which implies that $\mathrm{lct}(X,G)\leqslant 2$.

The surface $X$ can be obtained as a blow up $\pi\colon X\to
\mathbb{P}^{2}$ of the four points
$$
P_{1}=\big(1:-1:-1\big),\ P_{2}=\big(-1:1:-1\big),\ P_{3}=\big(-1:-1:1\big),\ P_{4}=\big(1:1:1\big),%
$$
of the plane $\mathbb{P}^{2}$. Let $W$ be the curve in
$\mathbb{P}^{2}$ that is given by the equation
$$
x^{6}+y^{6}+z^{6}+\big(x^{2}+y^{2}+z^{2}\big)\big(x^{4}+y^{4}+z^{4}\big)=12x^{2}y^{2}z^{2}\subset\mathbb{P}^{2}\cong\mathrm{Proj}\Big(\mathbb{C}[x,y,z]\Big),
$$
and $Z$ be its proper transform on $X$. Then $Z$ is
$\mathrm{S}_{5}$-invariant (see \cite{InKa05}) and $Z\sim
-2K_{X}$.

The curves $Z$ and $\sum_{i=1}^{10}L_{i}$ are the only
$\mathrm{S}_{5}$-invariant curves in $|-2K_{X}|$.

Let $\mathcal{P}$ be the pencil generated by $Z$ and
$\sum_{i=1}^{10}L_{i}$. It follows from \cite{Edge81a} that
$\mathcal{P}$ is $\mathrm{A}_{5}$-invariant, and there are exactly
$5$ singular curves in $\mathcal{P}$, which can be described in
the following way:
\begin{itemize}
\item the curve $\sum_{i=1}^{10}L_{i}$;%
\item two irreducible rational curves $R_{1}$ and $R_{2}$ that have $6$ nodes;%
\item two fibers $F_{1}$ and $F_{2}$ each consisting of $5$ smooth rational curves.%
\end{itemize}

We have $m=2$ and $k=6$ by \cite{RaSl02}. The smallest $G$-orbit
are $\mathrm{Sing}(R_{1})$ and $\mathrm{Sing}(R_{2})$ (see
\cite{InKa05}).

Suppose that $\mathrm{lct}(X,G)<2$. Then there is an effective
$G$-invariant $\mathbb{Q}$-divisor $D$ on the~quintic surface $X$
such that $\mathrm{LCS}(X,\lambda D)\ne\varnothing$  and
$D\equiv-K_{X}$, where $0<\lambda\in\mathbb{Q}$ such that
$\lambda<2$.

We may assume that the support of $D$ does not contain $R_{1}$ and
$R_{2}$ due to Remark~\ref{remark:convexity}, because both log
pairs $(X, R_{1})$ and $(X, R_{2})$ are log canonical. Now arguing
as in the proof of Lemma~\ref{lemma:m-over-r}, we see that either
$\mathrm{LCS}(X,\lambda D)=\mathrm{Sing}(R_{1})$ or
$\mathrm{LCS}(X,\lambda D)=\mathrm{Sing}(R_{2})$.

Without loss of generality we may assume that the locus $\mathrm{LCS}(X,\lambda
D)$ consists of the singular points of the curve $R_{1}$. Denote them as
$O_{1},\ldots,O_{6}$. Then  $\mathrm{mult}_{O_{i}}(D)\leqslant 5/6$, because
$$
10=R_{1}\cdot D\geqslant\sum_{i=1}^{6}\mathrm{mult}_{O_{i}}\big(D\big)\mathrm{mult}_{O_{i}}\big(R_{1}\big)\geqslant 12\mathrm{mult}_{O_{i}}\big(D\big).%
$$

Let $\pi\colon U\to X$ be a blow up of the points $O_{1},\ldots,O_{6}$.~Then
$$
K_{U}+2\bar{D}+\sum_{i=1}^{6}\Big(2\mathrm{mult}_{O_{i}}\big(D\big)-1\Big)E_{i}\equiv \pi^{*}\Big(K_{X}+2D\Big),%
$$
where $E_{i}$ is the $\pi$-exceptional curve such that $\pi(E_{i})=O_{i}$, and
$\bar{D}$ is the proper transform of~$D$~on the surface $U$. Then
$\mathrm{mult}_{Q_{i}}(\bar{D})>1-\mathrm{mult}_{O_{i}}(D)$ for some point
$Q_{i}\in E_{i}$, where $i=1,\ldots,6$.

Let $\Sigma$ be the $G$-orbit of the point $Q_{i}$. Then
$|\Sigma\cap E_{i}|\geqslant 2$, because the stabilizer of $O_{i}$
acts faithfully on its tangent space. We have $|\Sigma\cap
E_{i}|=2$, because $\mathrm{mult}_{O_{i}}(D)\leqslant 5/6$ and
$$
\mathrm{mult}_{O_{i}}\big(D\big)=E_{i}\cdot\bar{D}>\big|\Sigma\cap E_{i}\big|\Big(1-\mathrm{mult}_{O_{i}}\big(D\big)\Big).%
$$

Let $\bar{R}_{1}$ be the proper transform of the curve $R_{1}$ on the surface
$U$. Then
$$
\Sigma=\bar{R}_{1}\bigcap\Big(E_{1}\cup E_{2}\cup E_{3}\cup E_{4}\cup E_{5}\cup E_{5}\Big),%
$$
because the orbit of length $2$ of the action on $E_{i}$ of the stabilizer of
$O_{i}$ is unique. We have
$$
12\Big(1-\mathrm{mult}_{O_{i}}\big(D\big)\Big)=10-2\sum_{i=1}^{6}\mathrm{mult}_{O_{i}}\big(D\big)=\bar{R}_{1}\cdot\bar{D}\geqslant 2\left(\sum_{i=1}^{6}\mathrm{mult}_{Q_{i}}\big(\bar{D}\big)\right)>12\Big(1-\mathrm{mult}_{O_{i}}\big(D\big)\Big),%
$$
which is a contradiction.
\end{proof}

\begin{lemma}
\label{lemma:quintic-cyclic} Suppose that $K_{X}^{2}=5$ and
$G=\mathbb{Z}_{5}$. Then $\mathrm{lct}(X,G)=4/5$ holds.
\end{lemma}

\begin{proof}
It is well known that the group $G$ fixes exactly two points of the surfaces
$X$ (see \cite{RaSl02}), which we denote as $O_{1}$ and $O_{2}$. There are five
conics $Z_{1},\ldots,Z_{5}\subset X$~that passes through $O_{1}$, and the
divisor $\sum_{i=1}^{5}Z_{i}\sim -2K_{X}$ is $G$-invariant, which implies that
$\mathrm{lct}(X,G)\leqslant 4/5$.

Suppose that $\mathrm{lct}(X,G)<4/5$. Then there is an effective $G$-invariant
$\mathbb{Q}$-divisor $D$ on the~quintic surface $X$ such that
$\mathrm{LCS}(X,\lambda D)\ne\varnothing$  and $D\equiv-K_{X}$, where
$0<\lambda\in\mathbb{Q}$ such that $\lambda<4/5$.

The proof of Lemma~\ref{lemma:m-over-r} implies that $\mathrm{LCS}(X,\lambda
D)=\{O_{1}\}$ or $\mathrm{LCS}(X,\lambda D)=\{O_{1}\}$.

Without loss of generality, we may assume that $\mathrm{LCS}(X,\lambda
D)=\{O_{1}\}$, and we may assume that the support of the divisor $D$ does not
contain the conics $Z_{1},\ldots,Z_{5}$ by Remark~\ref{remark:convexity}. Then
$$
2=Z_{1}\cdot D\geqslant\mathrm{mult}_{O_{1}}\big(D\big).
$$

Let $\pi\colon U\to X$ be a blow up of the point $O_{1}$. and $E$
be the $\pi$-exceptional curve. Then
$$
\mathrm{mult}_{Q}\big(\bar{D}\big)\geqslant 2\big\slash\lambda-\mathrm{mult}_{O_{1}}\big(D\big)>5/2-\mathrm{mult}_{O_{1}}\big(D\big)%
$$
for some point $Q\in E$ by
Corollary~\ref{corollary:blow-up-inequality}, where $\bar{D}$ is
the proper transform of $D$ on the surface $U$.

The point $Q$ must be $G$-invariant, because otherwise
$$
\mathrm{mult}_{O_{1}}\big(D\big)=E\cdot\bar{D}>5\Big(5/2-\mathrm{mult}_{O_{1}}\big(D\big)\Big),
$$
which is impossible, because $\mathrm{mult}_{O_{1}}(D)\leqslant
2$.

Let $\bar{Z}_{i}$ be the proper transform of the conic $Z_{i}$ on
the surface $U$. Then $Q\not\in\cup_{i=1}^{5}\bar{Z}_{i}$, and
there is a birational morphism $\phi\colon U\to\mathbb{P}^{2}$
that contracts the curves $\bar{Z}_{1},\ldots,\bar{Z}_{5}$.

The curve $\phi(E)$ is a conic that contains
$\phi(\bar{Z}_{1}),\ldots,\phi(\bar{Z}_{5})$. Let $T_{i}$ be the
proper transform on the surface $U$ of the line in
$\mathbb{P}^{2}$ that passes through the points $\phi(Q)$ and
$\phi(\bar{Z}_{i})$. The log pair
$$
\left(X,\ \frac{\lambda}{3}\sum_{i=1}^{5}\pi(T_{i})\right)
$$
has log terminal singularities, and $\sum_{i=1}^{5}\pi(T_{i})\equiv 3D$. Thus,
we may assume that the support of the divisor $\bar{D}$ does not contain any of
the curves $T_{1},\cdots,T_{5}$ due to Remark~\ref{remark:convexity}. Then
$$
3-\mathrm{mult}_{O_{1}}\big(D\big)\geqslant T_{i}\cdot\bar{D}\geqslant\mathrm{mult}_{Q}\big(\bar{D}\big),%
$$
which implies that
$\mathrm{mult}_{O_{1}}(D)+\mathrm{mult}_{Q}(\bar{D})\leqslant 3$.

Let $\xi\colon V\to U$ be a blow up of the point $Q$, and $F$ be
the $\xi$-exceptional divisor. Then
$$
K_{W}+\lambda\grave{D}+\Big(\lambda\mathrm{mult}_{O_{1}}\big(D\big)-1\Big)\grave{E}+\Big(\lambda\mathrm{mult}_{O_{1}}\big(D\big)+\lambda\mathrm{mult}_{Q}\big(\bar{D}\big)-2\Big)F \equiv (\pi\circ\xi)^{*}\Big(K_{X}+\lambda D\Big),%
$$
where $\grave{D}$ and $\grave{E}$ are proper transforms of $D$ and $E$ on the
surface $V$, respectively. The log pair
$$
\Big(W,\ \lambda\grave{D}+\Big(\lambda\mathrm{mult}_{O_{1}}\big(D\big)-1\Big)\grave{E}+\Big(\lambda\mathrm{mult}_{O_{1}}\big(D\big)+\lambda\mathrm{mult}_{Q}\big(\bar{D}\big)-2\Big)F\Big)%
$$
is not log terminal at some point $P\in F$ by
Remark~\ref{remark:log-pull-back}, because $\mathrm{mult}_{O_{1}}(D)\leqslant
2$.

Suppose that $P\in\grave{E}$. Let $\grave{T}$ be the proper
transform on $V$ of the line on $\mathbb{P}^{2}$ that is tangent
to the conic $\phi(E)$ at the point $\phi(Q)$. Then
$P\in\grave{T}$, which implies that
$$
5-2\mathrm{mult}_{O_{1}}\big(D\big)-\mathrm{mult}_{Q}\big(\bar{D}\big)=\grave{T}\cdot\grave{D}\geqslant\mathrm{mult}_{P}\big(\grave{D}\big)>5-2\mathrm{mult}_{O_{1}}\big(D\big)-\mathrm{mult}_{Q}\big(\bar{D}\big),%
$$
because we may assume that $\grave{T}\not\subseteq\mathrm{Supp}(\grave{D})$ by
Remark~\ref{remark:convexity}. Hence, we have $P\not\in\grave{E}$.

The log pair
$(W,\lambda\grave{D}+(\lambda\mathrm{mult}_{O_{1}}(D)+\lambda\mathrm{mult}_{Q}(\bar{D})-2)F)$
is not log terminal at $P$. But
$$
\lambda\grave{D}+\Big(\lambda\mathrm{mult}_{O_{1}}\big(D\big)+\lambda\mathrm{mult}_{Q}\big(\bar{D}\big)-2\Big)F%
$$
is an effective divisor, because $\mathrm{mult}_{Q}(\bar{D})\geqslant
2/\lambda-\mathrm{mult}_{O_{1}}(D)$. Then
$$
\mathrm{mult}_{P}\big(\grave{D}\big)\geqslant
3\big\slash\lambda-\mathrm{mult}_{O_{1}}\big(D\big)-\mathrm{mult}_{Q}(\bar{D})>15/4-\mathrm{mult}_{O_{1}}\big(D\big)-\mathrm{mult}_{Q}(\bar{D}).
$$

Let $\grave{T}_{i}$ be the proper transform of $T_{i}$ on the
surface $V$. Suppose that $P\in\grave{T}_{k}$. Then
$$
3-\mathrm{mult}_{O_{1}}(D)-\mathrm{mult}_{Q}(\bar{D})=\grave{T}_{k}\cdot\grave{D}>15/4-\mathrm{mult}_{O_{1}}\big(D\big)-\mathrm{mult}_{Q}\big(\bar{D}\big),
$$
which is a contradiction. Thus, we see that
$P\not\in\cup_{i=1}^{5}\grave{T}_{i}$.

Let $M$ be an irreducible curve on $V$ such that $P\in M$, the
curve $\phi\circ\xi(M)$ is a line that passes  through the point
$\phi(Q)$. Then $\pi\circ\xi(M)$ has an ordinary double point at
$O_{1}$, and $\pi\circ\xi(M)\equiv -K_{X}$, because
$P\not\in\cup_{i=1}^{5}\grave{T}_{i}$. We may assume that
$M\not\subseteq\mathrm{Supp}(\grave{D})$ by
Remark~\ref{remark:convexity}. Then
$$
5-2\mathrm{mult}_{O_{1}}\big(D\big)-\mathrm{mult}_{Q}\big(\bar{D}\big)=M\cdot\grave{D}>15/4-\mathrm{mult}_{O_{1}}\big(D\big)-\mathrm{mult}_{Q}\big(\bar{D}\big),
$$
which implies that $\mathrm{mult}_{O_{1}}(D)\leqslant 5/4$. But
$\mathrm{mult}_{O_{1}}(D)>5/4$.
\end{proof}

We did not prove that groups in Example~\ref{example:quintic} and
Lemmaa~\ref{lemma:Fermat-cubic}, \ref{lemma:quintic-A5} and
\ref{lemma:quintic-cyclic} act on $X$ in such a way that
the~$G$-in\-va\-ri\-ant subgroup in $\mathrm{Pic}(X)$ is
$\mathbb{Z}$. But the latter is well--known (see \cite{DoIs06}).

\section{Direct products}
\label{section:products}

Let $X$ be an arbitrary smooth Fano variety, and let $G$ be a
finite subgroup in $\mathrm{Aut}(X)$ such that the $G$-invariant
subgroup of the group $\mathrm{Pic}(X)$ is $\mathbb{Z}$.

\begin{definition}
\label{definition:birational-rigidity} The variety $X$ is said to be
$G$-birationally superrigid if for every $G$-invariant linear system
$\mathcal{M}$ on the variety $X$ that does not have any fixed components,
the~singularities of the~log pair $(X,\lambda\mathcal{M})$~are~canonical, where
$\lambda\in\mathbb{Q}$ such~that $\lambda>0$ and
$K_{X}+\lambda\mathcal{M}\equiv 0$.
\end{definition}

The following result is well--known (see \cite{Ma67},
\cite{DoIs06}).

\begin{lemma}
\label{lemma:k-K-square} Suppose that $X$ is a smooth del Pezzo surface such
that
$$
\big|\Sigma\big|\geqslant K_{X}^{2}
$$
for any $G$-orbit $\Sigma\subset X$. Then $X$ is
$G$-birationally superrigid.
\end{lemma}

\begin{proof}
Suppose that the surface $X$ is not $G$-birationally superrigid. Then there is
a $G$-invariant linear system $\mathcal{M}$ on the surface $X$ such that
$\mathcal{M}$ does not have fixed curves, but $(X, \lambda\mathcal{M})$~is~not
canonical at some point $O\in X$, where $\lambda\in\mathbb{Q}$ such~that
$\lambda>0$ and $K_{X}+\lambda\mathcal{M}\equiv 0$.

Let $\Sigma$ be the $G$-orbit of the point $O$. Then
$\mathrm{mult}_{P}(\mathcal{M})>1/\lambda$ for every point $P\in \Sigma$. Then
$$
K_{X}^{2}\big\slash\lambda^{2}=M_{1}\cdot M_{2}\geqslant\sum_{P\in\Sigma}\mathrm{mult}^{2}_{P}\big(\mathcal{M}\big)>\big|\Sigma\big|\big\slash\lambda^{2}\geqslant K_{X}^{2}\big\slash\lambda^{2},%
$$
where $M_{1}$ and $M_{2}$ are sufficiently general curves in $\mathcal{M}$.
\end{proof}

\begin{example}
\label{example:quintic-A5} Let $X$ be a smooth del Pezzo surface
such that $K_{X}^{2}=5$. Then
$\mathrm{Aut}(X)\cong\mathrm{S}_{5}$, and the proof of
Lemma~\ref{lemma:quintic-A5} implies that the surface $X$ is
$\mathrm{A}_{5}$-birationally superrigid by
Lemma~\ref{lemma:k-K-square}.
\end{example}

Let $X_{i}$ be a smooth $G_{i}$-birationally superrigid Fano variety, where
$G_{i}$ is a an arbitrary finite subgroup of $\mathrm{Aut}(X_{i})$ such that
the~$G_{i}$-invariant subgroup of $\mathrm{Pic}(X_{i})$ is $\mathbb{Z}$, and
$i=1,\ldots,r$.

\begin{theorem}
\label{theorem:Pukhlikov} Suppose that $\mathrm{lct}(X_{i},
G_{i})\geqslant 1$ for every $i=1,\ldots,r$. Then
\begin{itemize}
\item there is no $G_{1}\times\cdots\times G_{r}$-equivariant birational map $\rho\colon X_{1}\times\cdots\times X_{r}\dasharrow\mathbb{P}^{n}$;%
\item every $G_{1}\times\cdots\times G_{r}$-equivariant birational automorphism of $X_{1}\times\cdots\times X_{r}$ is biregular;%
\item for any $G_{1}\times\cdots\times G_{r}$-equivariant dominant
map $\rho\colon X_{1}\times\cdots\times X_{r}\dasharrow Y$, whose
general fiber is rationally connected, there a commutative diagram
$$
\xymatrix{
X_{1}\times\cdots\times X_{r}\ar@{->}[d]_{\pi}\ar@{-->}[rrrrd]^{\rho}\\
X_{i_{1}}\times\cdots\times X_{i_{k}}\ar@{-->}[rrrr]_{\xi}&&&&Y}%
$$
where $\xi$ is a birational~map, $\pi$ is a natural projection,
and $\{i_{1},\ldots,i_{k}\}\subseteq\{1,\ldots,r\}$.
\end{itemize}
\end{theorem}

\begin{proof}
The required assertion follows from the proof of Theorem~1 in
\cite{Pu04d}.
\end{proof}

\begin{example}
\label{example:Valentiner} The simple group $\mathrm{A}_{6}$ is a group of
automorphisms of the sextic
$$
10x^{3}y^{3}+9zx^{5}+9zy^{5}+27z^{6}=45x^{2}y^{2}z^{2}+135xyz^{4}\subset\mathbb{P}^{2}\cong\mathrm{Proj}\Big(\mathbb{C}[x,y,z]\Big)
$$
which induces an embedding
$\mathrm{A}_{6}\subset\mathrm{Aut}(\mathbb{P}^{2})$ such that
$\mathrm{lct}(\mathbb{P}^{2},\mathrm{A}_{6})=2$
by~Lemma~\ref{lemma:m-over-r} (see \cite{Cra99}), and
$\mathrm{A}_{6}\times \mathrm{A}_{6}$ acts naturally on
$\mathbb{P}^{2}\times\mathbb{P}^{2}$. There is an induced
embedding $\mathrm{A}_{6}\times
\mathrm{A}_{6}\cong\Omega\subset\mathrm{Bir}(\mathbb{P}^{4})$
such that $\Omega$ is not conjugated to a subgroup of
$\mathrm{Aut}(\mathbb{P}^{4})$ by Lemma~\ref{lemma:k-K-square} and
Theorem~\ref{theorem:Pukhlikov}.
\end{example}

\end{document}